\documentclass[sn-mathphys]{sn-jnl} 

\jyear{2022}%

\theoremstyle{thmstyleone}%
\newtheorem{thm}{Theorem}
\newtheorem{prop}[thm]{Proposition}
\newtheorem{conj}[thm]{Conjecture}
\newtheorem{cor}[thm]{Corollary}

\theoremstyle{thmstylethree}
\newtheorem{defn}[thm]{Definition}
\newtheorem{nota}[thm]{Notation}

\begin{document}

\title[Spaces of flattenings of spheres]{Spaces of flattenings of spheres}


\author{\fnm{Olakunle S} \sur{Abawonse}}\email{o.abawonse@northeastern.edu}

\affil{\orgdiv{Department of Mathematics}, \orgname{Northeastern University}, \\ \orgaddress{\street{360 Huntington Ave}, \city{Boston}, \postcode{02115}, \state{Massachusetts}, \country{USA}}}


\abstract{The spaces of flattenings of a simplicial sphere played a key role in the study of existence and uniqueness of differentiable structures on a simplicial sphere. In this paper, we will establish that the spaces of flattenings of some simplicial spheres and show that they have the homotopy type of the orthogonal group. }




\maketitle

\section{Introduction}\label{sec1}

The space of flattenings $F(L)$ associated to a simplicial $k$-sphere $L$ is the space of all simplicial embeddings of the cone $cL$ of $L$ in $\mathbb{R}^{k+1}$. The space of flateenings of simplicial spheres was first studied by S.S Cairns \cite{cairns} in the early 1940s. In \cite{cairns}, Cairns proved that the space of all flattenings of a simplicial 2-sphere which have an orientation preserving isomorphism onto a given triangulation is path connected.

The space of flattenings associated to a simplicial 1-sphere  is homotopy equivalent to the orthogonal group $O(2)$ as we will prove in Theorem \ref{my_flat} \cite{kayman}. The only known positive result about the topology of $F(L)$ when $\dim(L) > 2$ is the theorem of N.H Kuiper \cite{kuiper:nh} where he proved that $F(L)$ has the homotopy type of $O(n+1)$ when $L$ is the boundary of the $(n+1)-$simplex. In \cite{milin}, Milin showed that there exist a simplicial sphere $L$ of dimension $3$ whose subset of $F(L)$ consisting of flattenings which have an orientation preserving isomorphism onto a given triangulation is not path connected.

Let $\Delta^n$ denote an $n$-simplex. We will establish that the space of flattenings associated to $\partial\Delta^n \ast \partial\Delta^1$ has the homotopy type of an orthogonal group. 

Let $L$ be a simplicial k-sphere and CF(L) denote the quotient of $F(L)$ by $\mathrm{GL}_{k+1}$ the invertible $(k+1)\times (k+1)$ matrices. To establish the above statement, we will first prove that $CF(L)$ is contractible. In Section \ref{or:flat}, we will introduce a poset $P(L)$ called \textit{The poset of oriented matroid flattenings} of $L$ and a poset stratification map $\pi: CF(L) \rightarrow P(L)$. We will prove that $P(L)$ is contractible when $L$ is either a simplicial 1-sphere, $\partial\Delta^{n+1}$ the boundary of an $(n+1)$-simplex or $\partial \Delta^1 \ast \partial\Delta^{n+1}$

In Section \ref{top:flat}, we will prove for the above simplicial spheres that $\{\overline{\pi^{-1}(M)} : M \in P(L)\}$ is a totally normal cellular decomposition of $CF(L)$. Theorem \ref{thm:total} will then conclude that $\|P(L)\|$ can be embedded in $CF(L)$ as a deformation retract.

\section{Oriented Matroids} \label{sec:back}

We will view elements of $\mathbb{R}^n$ as  $1\times n$ row vectors so that $X$ is the rowspace of a $r \times n$ matrix. Suppose $X \in \mathrm{Gr}(r, \mathbb{R}^n)$ so that $X = \mathrm{Rowspace}(v_1\; v_2 \; v_3 \ldots \; v_n)$. We consider the following function $\chi : [n]^r \rightarrow \{+, - , 0\}$  associated to $X$
	$$\chi(i_1, i_2, \ldots , i_r) = \mathrm{sign}(\det(v_{i_1} \; v_{i_2} \cdots v_{i_r}) )$$
The collection $\{\pm \chi\}$ is independent of the choice of basis vectors for $X$. The resulting functions $(\pm \chi)$ defines a rank $r$ \textit{oriented matroid}.

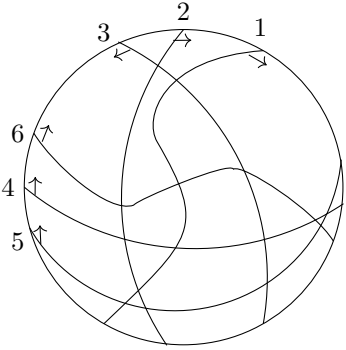
\begin{figure}[htb]
\begin{tikzpicture}
	\begin{scope}[scale=.7]
		\draw (0,0) circle(3);
		\draw[->](64:2.8)--(56:2.8);
		\draw[->](111:2.8)--(118:2.8);
		\draw[->, rotate=32](62:2.8)--(55:2.8);
		\draw[->, rotate=121](62:2.8)--(55:2.8);
		\draw[->, rotate=140](62:2.9)--(55:2.8);
		\draw[rotate=120] (60:3)  arc[radius = 5, start angle= 110, end angle= 184] node[above] at (0:3){$3$};
		\draw[rotate=240] (60:3)  arc[radius = 5, start angle= 110, end angle= 184] node[left] at (300:3){$4$};
		\draw[rotate=30] (60:3)  arc[radius = 5, start angle= 110, end angle= 184] node[above] at (60:3){$2$};
		\draw(195:3) arc[radius=3.15, start angle=210, end angle =355] node[left] at (200:3){$5$};
		
		\draw (60:3)..controls (100:2.5) and (120:1.5)..(120:1)..controls(-60:.9) and (-80:1)..(240:3) node[left] at (60:3.5){$1$};
\draw[] (160:3)..controls (170:2.5) and (200:1.5)..(200:1)..controls(195:.9) and (30:1)..(20:1)..controls(22:1.2) and (-10:2.5)..(-20:3) node[left] at (160:3){$6$};
\draw[->] (162:2.8)--(155:2.8);
\end{scope}
\end{tikzpicture}
\caption{Arrangement of Pseudospheres}
\label{pseudo}
\end{figure}

In general, an oriented matroid can be obtained from an arrangement of pseudospheres as evident by the following theorem. Figure \ref{pseudo} illustrates an arrangement of pseudospheres.

\begin{thm}(\cite{jim:law}){The Topological Representation Theorem (Folkman-Lawrence 1978)}
	The rank $r$ oriented matroids are exactly the sets $(E,\mathcal{V}^*)$ arising from essential {\em arrangements of pseudospheres} in $S^{r-1}$.
	 \end{thm}

A detailed introduction to the theory of oriented matroids can be found in the in the book \cite{anders:bjo}. Associated to a rank $r$ oriented matroid on $n$ elements are the functions $\pm \chi : [n]^r \rightarrow \{+, - , 0\}$ called the chirotopes. Let $\{+, -, 0\}$ be a poset with the partial order $0 < -$ and $0< +$.

\begin{defn}(\cite{macp:rob})
    Let $\mathcal{N} = (\pm \chi_1)$ and $\mathcal{M} = (\pm \chi_2)$ be two rank $r$ oriented matroids. We say that $\mathcal{N} \leq \mathcal{M}$ if and only if $\chi_1 \leq \chi_2$ or $\chi_1 \leq -\chi_2$. The oriented matroid $\mathcal{M}$ is said to {\em weak map} to $\mathcal{N}$. 
\end{defn}

 \begin{defn}(\cite{macp:rob})
 $\mathrm{MacP}(p,n)$ denotes the poset of all rank $p$ oriented matroids on elements $\{1,2,\ldots, n\}$, with weak map as the partial order. The poset is called the \textit{MacPhersonian} ~\cite{macp:rob}.
 \end{defn}

We have explained how to obtain a rank $r$ oriented matroid on $n$ elements from a rank $r$ subspace of $\mathbb{R}^n$. That is, there is  a function $\mu: \mathrm{Gr}(r, \mathbb{R}^n) \rightarrow \mathrm{MacP}(r,n): X \to (\pm \chi_X)$.

The following Proposition and Theorem are from the work of the author in \cite{kay:ab}, \cite{kayman}.

	\begin{prop}(\cite{kay:ab})\label{ref1}
		Let $M \in \mbox{MacP}(2,n)$. Then $\partial \overline{\mu^{-1}(M)} = \bigcup_{N < M} \mu^{-1}(N)$
	\end{prop}

\begin{thm}(\cite{kay:ab})\label{ref2} $\{\overline{\mu^{-1}(M)} : M \in \; \mbox{MacP}(2, n)\}$ is a regular cell decomposition of $Gr(2,\mathbb{R}^n)$.
	\end{thm}

\section{Cellular stratified spaces }\label{cellular}
 
\begin{defn}(\cite{Dai:Tam})
     A \textit{globular} $n$-cell is a subset $D$ of $D^n$ containing $H= \mathrm{Int}(D^n)$. We call $D \cap \partial D^n$ the \textit{boundary} of $D$ and denote it by $\partial D$. The number $n$ is called the \textit{globular dimension} of $D$. 
\end{defn}

A globular $n$-cell was introduced by Tamaki \cite{Dai:Tam} as an extension of closure of $n$-cells to non-closed cells.

\begin{defn} (\cite{Dai:Tam})
    Let $X$ be a topological space. For a non-negative integer $n$, an $n$-cell structure on a subspace $e \subset X$ is a pair $(D, \varphi)$ of a globular $n$-cell $D$ and a continuous map $$\varphi: D \rightarrow X$$
    satsifying the following conditions:
    \begin{itemize}
        \item $\varphi(D) = \bar{e}$ and $\varphi : D \rightarrow \bar{e}$ is a quotient map.
        
        \item The restriction $\varphi$ : $H \rightarrow e$ is a homeomorphism.
    \end{itemize}
    
\end{defn}

\begin{defn}(\cite{Dai:Tam})
Let $X$ be a topological space and $P$ be a poset with the Alexandroff topology. A stratification of $X$ indexed by $P$ is an open continuous map 
$$\pi : X \rightarrow P$$
satisfying the condition that for each $\lambda \in P$, $e_\lambda = \pi^{-1}(\lambda)$ is connected and locally closed. $X$ is called a \textit{cellular stratified space} if each $e_\lambda$ is homeomorphic to an open ball.
\end{defn}

\begin{defn}(\cite{Dai:Tam}, \cite{Fur:Muk})
    Let $X$ be a cellular stratifed space. $X$ is called totally normal if for each globular $n$-cell $(D_\lambda, \varphi)$,  and $e_\lambda = \varphi(\mathrm{Int}(D_\lambda))$
    \begin{enumerate}[(i)]
        \item If $e_\lambda \cap \overline{e_\mu} \neq \emptyset $, then $e_\lambda \subseteq \overline{e_\mu}$. 
        \item There exists a structure of a regular cell complex on $S^{n-1}$ containing $\partial D_\lambda$ as a cellular stratified subspace of $S^{n-1}$.
        \item For each cell $e$ in the cellular stratification on $\partial D_\lambda$, there exists a cell $e_\eta$ in $X$  and a map $b: D_\eta \rightarrow \partial D_\lambda$ such that $b(\mathrm{Int}(D_\eta)) = e$ and $\varphi_\lambda \circ b = \varphi_\eta$.
    \end{enumerate}
\end{defn}

\begin{thm}\label{thm:total} (\cite{Dai:Tam})
For a totally  normal cellular stratified space $X$ with stratification $\pi : X \rightarrow P$, there is an embedding of $\|P\|$ as a strong deformation retract of $X$. 
\end{thm}

\section{Flattenings} \label{def:flattenings}

\begin{defn} (\cite{losik:1}, \cite{milin})
	Let $L$ be a triangulation of a $k$-sphere, and let $cL$ be a simplicial cone over $L$. A flattening of $L$ is an embedding $\psi : cL \rightarrow \mathbb{R}^{k+1}$ that maps the cone vertex to the origin and it is linear on simplices of $cL$.  
\end{defn}

\begin{nota}
Let $L$ be a simplicial $k$-sphere. We denote as in \cite{losik:1} by $F(L)$ the space of all flattenings of $L$. Also, the group $GL_{k+1}$ of invertible $(k+1) \times (k+1)$ matrices acts on $F(L)$; the quotient space denoted by $CF(L)$ is the configuration space of $L$.
\end{nota}

The space $F(L)$ is an open subset of $\mathbb{R}^{(k+1)|\mbox{Vert}(L)|}$, and so has a natural smooth manifold structure. The space of flattenings comes up in the problem of existence and uniqueness of differentiable structures on triangulated manifolds (see \cite{losik:2}, \cite{kuiper:nh}). 

We will show that $CF(L)$ is contractible  when $L$ is a simplicial $1$-sphere, and so, $F(L)$ has the homotopy type of $O(2)$. Some few other non-trivial results that are known about the topology of $CF(L)$ and $F(L)$ are as follows.

\begin{thm}(\cite{cairns})
	Let $L$ be a triangulated $2$-sphere. Then $CF(L)$ is path connected. 
\end{thm}

For $\dim(L) \geq 3$, Cairns \cite{cairns} also showed that $CF(L)$ can be empty. When the dimension of $L$  is greater than $2$,  Milin \cite{milin} obtained the following negative result about the topology of $CF(L)$.
 
 \begin{thm}\cite{milin}
 	There exists a $3$ dimensional simplicial sphere $L$ such that $CF(L)$ is disconnected.
 \end{thm}
 
So far, for $n > 2$ the only known positive result about the homotopy type of $F(L)$ is the following result of Kuiper.

\begin{thm}(\cite{kuiper:nh})\label{kuip}
	Let $\partial \Delta^{n+1}$ be the boundary of an $(n+1)$-simplex. Then $F(\partial \Delta^{n+1})$ has the homotopy type of $O(n+1)$. 
\end{thm}
\begin{cor}(\cite{kuiper:nh})
	Let $\partial \Delta^{n+1}$ be the boundary of an $(n+1)$-simplex. Then any two smoothings of $\partial \Delta^{n+1}$ are diffeomorphic. 
\end{cor}

As in Theorem \ref{kuip}, we also obtain the following positive result for the simplicial sphere $\partial \Delta^1 \ast \partial \Delta^{n+1}$.

\begin{thm}\label{my_flat}
Let $\partial \Delta^{n+1}$ be the boundary of an $(n+1)$-simplex . Then $F(\partial \Delta^1 \ast \partial \Delta^{n+1})$ has the homotopy type of $O(n+2)$. Let $L$ be a simplicial $1$-sphere. Then $F(L)$ has the homotopy type of $O(2)$. 
\end{thm}

\begin{cor}\label{cor_flat}
	Let $\partial \Delta^{n+1}$ be the boundary of an $(n+1)$-simplex. Then any two smoothings of $\partial \Delta^1 \ast \partial \Delta^{n+1}$ are diffeomorphic. 
\end{cor}

\section{Oriented matroid flattenings} \label{or:flat}
Let $L$ be a triangulated of a $k$-sphere, and $\psi : cL \rightarrow \mathbb{R}^{k+1}$ a flattening of $L$. Then the arrangement of vectors $(\psi(v): v \in \mathrm{Vert}(L))$ determines a rank $k+1$ oriented matroid $M$. Definition \ref{comb_abstr} gives a combinatorial abstraction for oriented matroids obtained from flattenings of a simplicial sphere.

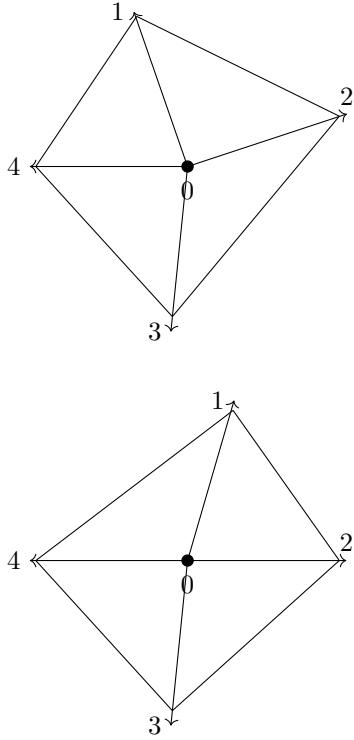
\begin{figure}[htb]
\centering
\begin{subfigure}[t]{0.35\textwidth}
    \begin{tikzpicture}[line join=bevel,z=-5.5]
    \coordinate (A1) at (-0.2,-2);
    \coordinate (A2) at (-2,0);
    \coordinate (A3) at (-0.67,2);
    \coordinate (A5) at (0,0);
    \coordinate (A4) at (2,0.67);
    \draw (A1) -- (A2) -- (A3) ;
    \draw (A1) -- (A4) -- (A3);
    \path[->] (0,0) edge node[at end, left]{$3$} (-0.22, -2.2);
    \path[->] (0,0) edge node[at end, left]{$4$} (-2.08, 0);
    \path[->] (0,0) edge node[at end, left]{$1$} (-0.707, 2.06);
    \path[->] (0,0) edge node[at end, above]{$2$} (2.1, 0.7);
    \foreach \Coor/\Texto/\Pos in 
    {A5/0/below
    }
    \node[circle,draw,inner sep=1.5pt,fill=black,label={\Pos:$\Texto$}] 
    at (\Coor) {};
    \end{tikzpicture}
\end{subfigure}

~

\begin{subfigure}[t]{0.35\textwidth}
    \begin{tikzpicture}[line join=bevel,z=-5.5]
    \coordinate (A1) at (-0.2,-2);
    \coordinate (A2) at (-2,0);
    \coordinate (A3) at (0.6,2);
    \coordinate (A5) at (0,0);
    \coordinate (A4) at (2,0);
    \draw (A1) -- (A2) -- (A3) ;
    \draw (A1) -- (A4) -- (A3);
    \path[->] (0,0) edge node[at end, left]{$3$} (-0.22, -2.2);
    \path[->] (0,0) edge node[at end, left]{$4$} (-2.08, 0);
    \path[->] (0,0) edge node[at end, left]{$1$} (0.613, 2.133);
    \path[->] (0,0) edge node[at end, above]{$2$} (2.093, 0);
    \foreach \Coor/\Texto/\Pos in 
    {A5/0/below
    }
    \node[circle,draw,inner sep=1.5pt,fill=black,label={\Pos:$\Texto$}] 
    at (\Coor) {};
    \end{tikzpicture}
\end{subfigure}
\caption{Flattenings of a simplicial $1$-sphere.}
\end{figure}

\begin{defn}\label{comb_abstr}
    Let $L$ be a simplicial sphere of dimension $k$. An oriented matroid flattening of $L$ is a rank $k+1$ oriented matroid $\mathcal{M}$ satisfying the following:
    \begin{enumerate}[(i)]
        \item The elements of $\mathcal{M}$ are the vertices of $L$.
        \item The set of vertices in a simplex are independent.
        \item The set of vertices in a simplex has no other elements in its convex hull.
        
    \end{enumerate}
\end{defn}

 \begin{nota}
 The poset of all oriented matroid flattenings of $L$ is denoted by $P(L)$. 
\end{nota}

\begin{prop}\label{contr_prop}
Let $L$ be a simplicial sphere and $\partial \Delta^n$ the boundary of an $n$-simplex. Then $\|P(L)\|$ is contractible when $L$ is either a simplicial $1$-sphere, $\partial \Delta^n$ or $\partial \Delta^1 \ast \partial \Delta^n$. 
\end{prop}

\begin{proof}
The poset $P(\partial \Delta^n)$ consists of a point. Let $P(\partial \Delta^n) = \{\mathcal{M}_n\}$. For the sphere $\partial\Delta^1 \ast \partial \Delta^n$, $P(\partial\Delta^1 \ast \partial \Delta^n)$ has a minimum; given by the join $\mathcal{M}_1 \oplus \mathcal{M}_n$ of two oriented matroids.

In the case when $L$ is a simplicial $1$-sphere, this will follow by induction on the number of vertices in $\mathrm{Vert}(L)$. Let $L_n$ denote a simplical $1$-sphere on $n$ vertices. We know that $P(L_3)$ consists of a point say $M_0 = (\pm \chi_0)$. In the following argument, we will consider chirotopes with positive value on  the basis $\{1,2\}$ 

Let $\Sigma^{n+1}$ denote a subposet of $P(L_{n+1})$ consisting of $\mathcal{M}'$ such that $\mathcal{M}'\setminus \{n+1\}$ is an element of $P(L_n)$. An oriented matroid in $\Sigma^{n+1}$ is thus an extension of an oriented matroid $\mathcal{M}$ in $P(L_n)$ by an element ${n+1}$, with $n+1$ lying in the convex hull of $\{1,n\}$.

There is a poset map $P(L_{n+1}) \rightarrow \Sigma^{n+1}$ obtained as composition of some poset maps  as given below.
Let $f_0 : P(L_{n+1}) \rightarrow P(L_{n+1})$ defined as:
$$f_0(\chi)(B) = \left\{\begin{array}{ccc}
    \chi(B) & \mbox{if} & B \neq (n, 1)  \\
     0 &  \mbox{if} & B = (n,1) \; \mbox{and}\; \chi(n, 1) \in \{0, -\}\\
     + & \mbox{if} & B = (n,1) \; \mbox{and} \; \chi(n,1) = +
\end{array}\right\}$$
Let $P_0 = f_0(P(L_{n+1}))$. The poset map $f_0$ is a lowering homotopy, and so $\|P_0\|$ is homotopy equivalent to $\|P(L_{n+1})\|$. We again consider another poset map $f_1: P_0 \rightarrow P_0$ defined as:

$$f_1(\chi)(B) = \left\{\begin{array}{ccc}
    \chi(B) & \mbox{if} & B \neq (n,1)  \\
     + & \mbox{if} & B = (n,1) 
\end{array}\right\}$$

The image of $f_1$ is denoted is given by $f_1(P_0) = \Sigma^{n+1}$. The poset map $f_1 : P_0 \rightarrow P_0$ is a raising homotopy, and so $\|P_0\|$ is homotopy equivalent to $\|\Sigma^{n+1}\|$. The poset map $\Sigma^{n+1} \rightarrow P(L_n)$ induces a homotopy equivalence between $\|\Sigma^{n+1}\|$ and $\|P(L_n)\|$.
\end{proof}

For a simplicial sphere $L$, there is a stratification map $\mu_0: CF(L) \rightarrow P(L)$. 

\begin{conj}
	Let $L$ be a simplicial sphere of dimension at least $2$. Then $\|P(L)\|$ is contractible.
\end{conj}

\section{Topology of space of flattenings of some spheres} \label{top:flat}
Let $\mu' : \mathrm{Gr}(r, \mathbb{R}^{r+2}) \rightarrow \mathrm{MacP}(r, r+2)$ and $\mu: \mathrm{Gr}(2, \mathbb{R}^n) \rightarrow \mathrm{MacP}(2,n) $.  Let $\mu_0 : CF(L) \rightarrow P(L)$ be the restriction of $\mu'$ to $CF(L)$ when $L = \partial \Delta^1 \ast \partial \Delta^{r-1}$ or the restriction of $\mu$ when $L$ is a simplicial 1-sphere on $n$ vertices.

The stratification map $\mu_0: CF(L) \rightarrow P(L)$ gives a decomposition of $CF(L)$ into semi-algebraic sets $\{\mu_0^{-1}(M) : M \in P(L)\}$. When $L$ is a simplicial 1-sphere or $L = \partial \Delta^1 \ast \partial \Delta^n$, we will show that the decomposition is a totally normal cellular decomposition.

We have the following commutative diagram

\begin{tikzcd}
\mathrm{Gr}(r, \mathbb{R}^{r+2}) \arrow[r, "\mu'"] \arrow[d, "V \mapsto V^{\perp}", labels = left ] & \mathrm{MacP}(r, r+2) \arrow[d, "M\mapsto M^*"]\\
\mathrm{Gr}(2, \mathbb{R}^{r+2}) \arrow[r, "\mu"] & \mathrm{MacP}(2, r+2)
\end{tikzcd}

The commutativity of the diagram follows from the fact that  

$V = (I_r | A) \in \mathrm{Gr}(r, \mathbb{R}^{r+2})$ if and only if $V^{\perp} = \mathrm{Rowspace}(-A^T|I_2) \in \mathrm{Gr}(2, \mathbb{R}^{r+2})$. The oriented matroid $M^*$ is called the dual of $M$.

The map $\mathrm{Gr}(r, \mathbb{R}^{r+2}) \rightarrow \mathrm{Gr}(2, \mathbb{R}^{r+2}) : V \mapsto V^{\perp}$ is a homeomorphism and the poset map $\mathrm{MacP}(r, r+2) \rightarrow \mathrm{MacP}(2, r+2) : M \mapsto M^*$ is a poset isomorphism.

The  following result thus follows  from Theorem \ref{ref2}
and the  commutativity of the diagram described above.

\begin{thm}\label{thm:reg}
	Let $M \in \mbox{MacP}(r, r+2)$ be a rank $r$ oriented matroid on $r+2$ elements, and $\mu' : \mathrm{Gr}(r, \mathbb{R}^{r+2}) \rightarrow \mbox{MacP}(r,r+2) $. Then $\{\overline{(\mu')^{-1}(M)}: M \in \mathrm{MacP}(r, r+2)\}$ is a regular cell decomposition of $\mathrm{Gr}(r, \mathbb{R}^{r+2})$.
\end{thm}

\begin{prop}\label{tot_nor}
Let $L$ be a simplicial sphere and  $\mu_0 : CF(L) \rightarrow P(L)$ a stratification map. If $L$ is a simplicial $1$-sphere or $L = \partial \Delta^1 \ast \partial \Delta^n$, then the decomposition $\{\mu_0^{-1}(M) : M \in P(L)\}$ is a totally normal cellular decomposition of $CF(L).$
\end{prop}

\begin{proof}
\begin{enumerate}[(i)]
    \item Suppose $L$ is as given above. It was proven in Proposition \ref{ref1} that if $N, M \in P(L)$ such that $N < M$, then $\mu_0^{-1}(N) \subseteq \overline{\mu_0^{-1}(M)}$. So, the decomposition  $\{\mu_0^{-1}(M) : M \in P(L)\}$ is normal.
    
    \item In Theorem \ref{ref2}, it was proven that $\{\overline{\mu^{-1}(M)} : M \in \; \mbox{MacP}(2, |\mathrm{Vert}(L)|)\}$ is a regular cell decomposition of $Gr(2,\mathbb{R}^{|\mathrm{Vert}(L)|})$. Similarly, we have in Theorem \ref{thm:reg} that $\{\overline{(\mu')^{-1}(M)} : M \in \; \mbox{MacP}(r, r+2)\}$ is a regular cell decomposition of $Gr(r,\mathbb{R}^{r+2})$.
    
    If $L$ is a simplicial $1$-sphere, and $M \in P(L)$, then $\partial \overline{\mu^{-1}(M)}$ is a regular cellular cell complex  homeomorphic to a sphere. Let $\overline{\mu_0^{-1}(M)}$ denote the closure of $\mu_0^{-1}(M)$ in $CF(L)$. Then $\partial \overline{\mu^{-1}(M)}$ contains $\partial \overline{\mu_0^{-1}(M)}$ as a cellular stratified subspace. Similarly when $L = \partial \Delta^1 \ast \partial \Delta^n$ and $M \in P(L)$, $\partial \overline{(\mu')^{-1}(M)}$ contains $\partial \overline{\mu_0^{-1}(M)}$ as a cellular stratified subspace.
    
    \item $D_M = \overline{\mu_0^{-1}(M)}$, and let $\varphi_M$ be the restriction to $D_M$ of the characteristic map of the cell $\overline{\mu^{-1}(M)}$ if $L$ is a simplicial $1$-sphere or  restriction of the characteristic map of $\overline{(\mu')^{-1}(M)}$ if $L = \partial \Delta^1 \ast \partial \Delta^n$.
    
    For a cell $e$ in the boundary of $D_M$, there exists an oriented matroid $N$ in $P(L)$ such that $N < M$ and $\mu_0^{-1}(N) = e$. The map $b : D_N \rightarrow \partial D_M$ is given by $b = (\varphi_M)^{-1} \circ \varphi_N.$ 
    
    \end{enumerate}
\end{proof}

\begin{proof}[Proof of Theorem \ref{my_flat}]
Suppose $L$ is a simplicial $1$-sphere or $L= \partial \Delta^1 \ast \partial \Delta^n$. The decomposition $\{\mu_0^{-1}(M): M \in P(L)\}$ is a totally normal cellular decomposition of $CF(L)$ by Proposition \ref{tot_nor}. It thus follows from Theorem \ref{thm:total} that $\|P(L)\|$ is a deformation retract of $CF(L)$. We know from  Proposition \ref{contr_prop} that $\|P(L)\|$ is contractible. Hence, $CF(L)$ is contractible.

Suppose $L$ is a simplicial $1$-sphere. We know that $F(L)|_H \cong \mathrm{GL}_2(\mathbb{R}) \times CF(L)$. Hence $F(L)$ has the homotopy type of $O(2)$. Similarly, if $L = \partial \Delta^1 \ast \partial \Delta^n$, then $F(L) \cong \mathrm{GL}_{n+1}(\mathbb{R}) \times CF(L)$. Hence $F(L)$ has the homotopy type of $O(n+1)$.
\end{proof}

\bibliography{MacP}


\begin{thebibliography}{12}
\ifx \bisbn   \undefined \def \bisbn  #1{ISBN #1}\fi
\ifx \binits  \undefined \def \binits#1{#1}\fi
\ifx \bauthor  \undefined \def \bauthor#1{#1}\fi
\ifx \batitle  \undefined \def \batitle#1{#1}\fi
\ifx \bjtitle  \undefined \def \bjtitle#1{#1}\fi
\ifx \bvolume  \undefined \def \bvolume#1{\textbf{#1}}\fi
\ifx \byear  \undefined \def \byear#1{#1}\fi
\ifx \bissue  \undefined \def \bissue#1{#1}\fi
\ifx \bfpage  \undefined \def \bfpage#1{#1}\fi
\ifx \blpage  \undefined \def \blpage #1{#1}\fi
\ifx \burl  \undefined \def \burl#1{\textsf{#1}}\fi
\ifx \doiurl  \undefined \def \doiurl#1{\url{https://doi.org/#1}}\fi
\ifx \betal  \undefined \def \betal{\textit{et al.}}\fi
\ifx \binstitute  \undefined \def \binstitute#1{#1}\fi
\ifx \binstitutionaled  \undefined \def \binstitutionaled#1{#1}\fi
\ifx \bctitle  \undefined \def \bctitle#1{#1}\fi
\ifx \beditor  \undefined \def \beditor#1{#1}\fi
\ifx \bpublisher  \undefined \def \bpublisher#1{#1}\fi
\ifx \bbtitle  \undefined \def \bbtitle#1{#1}\fi
\ifx \bedition  \undefined \def \bedition#1{#1}\fi
\ifx \bseriesno  \undefined \def \bseriesno#1{#1}\fi
\ifx \blocation  \undefined \def \blocation#1{#1}\fi
\ifx \bsertitle  \undefined \def \bsertitle#1{#1}\fi
\ifx \bsnm \undefined \def \bsnm#1{#1}\fi
\ifx \bsuffix \undefined \def \bsuffix#1{#1}\fi
\ifx \bparticle \undefined \def \bparticle#1{#1}\fi
\ifx \barticle \undefined \def \barticle#1{#1}\fi
\bibcommenthead
\ifx \bconfdate \undefined \def \bconfdate #1{#1}\fi
\ifx \botherref \undefined \def \botherref #1{#1}\fi
\ifx \url \undefined \def \url#1{\textsf{#1}}\fi
\ifx \bchapter \undefined \def \bchapter#1{#1}\fi
\ifx \bbook \undefined \def \bbook#1{#1}\fi
\ifx \bcomment \undefined \def \bcomment#1{#1}\fi
\ifx \oauthor \undefined \def \oauthor#1{#1}\fi
\ifx \citeauthoryear \undefined \def \citeauthoryear#1{#1}\fi
\ifx \endbibitem  \undefined \def \endbibitem {}\fi
\ifx \bconflocation  \undefined \def \bconflocation#1{#1}\fi
\ifx \arxivurl  \undefined \def \arxivurl#1{\textsf{#1}}\fi
\csname PreBibitemsHook\endcsname

\bibitem{cairns}
\begin{barticle}
\bauthor{\bsnm{Cairns}, \binits{S.S.}}:
\batitle{Isotopic deformations of geodesic complexes on the 2-sphere and on the
  plane}.
\bjtitle{Ann. of Math. (2)}
\bvolume{45},
\bfpage{207}--\blpage{217}
(\byear{1944}).
\doiurl{10.2307/1969263}
\end{barticle}
\endbibitem

\bibitem{kayman}
\begin{botherref}
\oauthor{\bsnm{Abawonse}, \binits{O.S.}}:
On the topology of flags of oriented matroids and spaces of flattenings of
  spheres.
PhD thesis
(2022).
Copyright - Database copyright ProQuest LLC; ProQuest does not claim copyright
  in the individual underlying works; Last updated - 2022-08-19.
\url{https://link.ezproxy.neu.edu/login?url=https://www.proquest.com/dissertations-theses/on-topology-flags-oriented-matroids-spaces/docview/2699687391/se-2}
\end{botherref}
\endbibitem

\bibitem{kuiper:nh}
\begin{bchapter}
\bauthor{\bsnm{Kuiper}, \binits{N.H.}}:
\bctitle{On the smoothings of trangulated and combinatorial manifolds}.
In: \bbtitle{Differential and {C}ombinatorial {T}opology ({A} {S}ymposium in
  {H}onor of {M}arston {M}orse)},
pp. \bfpage{3}--\blpage{22}.
\bpublisher{Princeton Univ. Press, Princeton, N.J.}, \blocation{???}
(\byear{1965})
\end{bchapter}
\endbibitem

\bibitem{milin}
\begin{barticle}
\bauthor{\bsnm{Milin}, \binits{L.}}:
\batitle{A combinatorial computation of the first {P}ontryagin class of the
  complex projective plane}.
\bjtitle{Geom. Dedicata}
\bvolume{49}(\bissue{3}),
\bfpage{253}--\blpage{291}
(\byear{1994}).
\doiurl{10.1007/BF01264030}
\end{barticle}
\endbibitem

\bibitem{jim:law}
\begin{barticle}
\bauthor{\bsnm{Folkman}, \binits{J.}},
\bauthor{\bsnm{Lawrence}, \binits{J.}}:
\batitle{Oriented matroids}.
\bjtitle{J. Combin. Theory Ser. B}
\bvolume{25}(\bissue{2}),
\bfpage{199}--\blpage{236}
(\byear{1978}).
\doiurl{10.1016/0095-8956(78)90039-4}
\end{barticle}
\endbibitem

\bibitem{anders:bjo}
\begin{barticle}
\bauthor{\bsnm{Bj\"{o}rner}, \binits{A.}}:
\batitle{Posets, regular {CW} complexes and {B}ruhat order}.
\bjtitle{European J. Combin.}
\bvolume{5}(\bissue{1}),
\bfpage{7}--\blpage{16}
(\byear{1984}).
\doiurl{10.1016/S0195-6698(84)80012-8}
\end{barticle}
\endbibitem

\bibitem{macp:rob}
\begin{barticle}
\bauthor{\bsnm{MacPherson}, \binits{R.}}:
\batitle{Combinatorial differential manifolds, in: {"Topological Methods in
  Modern Mathematics: A Symposium in Honor of John Milnor's Sixtieth Birthday
  Stony Brook NY 1991"}}.
\bjtitle{Publish or Perish, Houston TX 1993}
\bvolume{(225282)},
\bfpage{203}--\blpage{221}
(\byear{1993}).
\doiurl{10.2307/1970177}
\end{barticle}
\endbibitem

\bibitem{kay:ab}
\begin{botherref}
\oauthor{\bsnm{Abawonse}, \binits{O.S.}}:
Homeomorphism type of combinatorial grassmannian and flag manifold
(2022).
Preprint at \url{https://arxiv.org/abs/2205.09553}
\end{botherref}
\endbibitem

\bibitem{Dai:Tam}
\begin{bchapter}
\bauthor{\bsnm{Tamaki}, \binits{D.}}:
\bctitle{Cellular stratified spaces}.
In: \bbtitle{Combinatorial and Toric Homotopy}.
\bsertitle{Lect. Notes Ser. Inst. Math. Sci. Natl. Univ. Singap.},
vol. \bseriesno{35},
pp. \bfpage{305}--\blpage{435}.
\bpublisher{World Sci. Publ., Hackensack, NJ}, \blocation{???}
(\byear{2018})
\end{bchapter}
\endbibitem

\bibitem{Fur:Muk}
\begin{barticle}
\bauthor{\bsnm{Furuse}, \binits{M.}},
\bauthor{\bsnm{Mukouyama}, \binits{T.}},
\bauthor{\bsnm{Tamaki}, \binits{D.}}:
\batitle{Totally normal cellular stratified spaces and applications to the
  configuration space of graphs}.
\bjtitle{Topol. Methods Nonlinear Anal.}
\bvolume{45}(\bissue{1}),
\bfpage{169}--\blpage{214}
(\byear{2015}).
\doiurl{10.12775/TMNA.2015.010}
\end{barticle}
\endbibitem

\bibitem{losik:1}
\begin{barticle}
\bauthor{\bsnm{Gabri\`elov}, \binits{A.M.}},
\bauthor{\bsnm{Gel\cprime~fand}, \binits{I.M.}},
\bauthor{\bsnm{Losik}, \binits{M.V.}}:
\batitle{A local combinatorial formula for the first {P}ontrjagin class}.
\bjtitle{Funkcional. Anal. i Prilo\v{z}en.}
\bvolume{10}(\bissue{1}),
\bfpage{14}--\blpage{17}
(\byear{1976})
\end{barticle}
\endbibitem

\bibitem{losik:2}
\begin{barticle}
\bauthor{\bsnm{Gabri\`elov}, \binits{A.M.}},
\bauthor{\bsnm{Gel\cprime~fand}, \binits{I.M.}},
\bauthor{\bsnm{Losik}, \binits{M.V.}}:
\batitle{Combinatorial computation of characteristic classes. {I}, {II}}.
\bjtitle{Funkcional. Anal. i Prilo\v{z}en.}
\bvolume{9}(\bissue{2}),
\bfpage{12}--\blpage{28919753526}
(\byear{1975})
\end{barticle}
\endbibitem

\end{thebibliography}


\end{document}